\documentclass{birkjour}
\usepackage{amssymb,bbm,hyperref,color}

\newtheorem{theorem}{Theorem}[section]
\newtheorem{lemma}[theorem]{Lemma}
\newtheorem{corollary}[theorem]{Corollary}

\theoremstyle{remark}

\numberwithin{equation}{section}

\title[The essential norms of Toeplitz operators]
{The essential norms of Toeplitz operators 
with symbols in \boldmath{$C+H^\infty$} on weighted Hardy spaces 
are independent of the weights}

\author[O. Karlovych]{Oleksiy Karlovych}
\address{
Centro de Matem\'atica e Aplica\c{c}\~oes\\
Departamento de Matem\'atica\\
Faculdade de Ci\^encias e Tecnologia\\
Universidade Nova de Lisboa\\
Quinta da Torre\\
2829--516 Caparica\\
Portugal} \email{oyk@fct.unl.pt}
\author[E. Shargorodsky]{Eugene Shargorodsky}
\address{%
Department of Mathematics\\
King's College London\\
Strand, London WC2R 2LS\\
United Kingdom}
\email{eugene.shargorodsky@kcl.ac.uk}
\begin{document}
\begin{abstract}
Let $1<p<\infty$, let $H^p$ be the Hardy space on the unit circle, and let 
$H^p(w)$ be the Hardy space with a Muckenhoupt weight $w\in A_p$ on the unit 
circle. In 1988, B\"ottcher, Krupnik and Silbermann proved that
the essential norm of the Toeplitz operator $T(a)$ with $a\in C$
on the weighted Hardy space $H^2(\varrho)$ with a power weight $\varrho\in A_2$
is equal to $\|a\|_{L^\infty}$. This implies that the essential norm of
$T(a)$ on $H^2(\varrho)$ does not depend on $\varrho$. We extend this result 
and show that if $a\in C+H^\infty$, then, for $1<p<\infty$, the essential 
norms of the Toeplitz operator $T(a)$ on $H^p$ and on $H^p(w)$ are the same 
for all $w\in A_p$. In particular, if $w\in A_2$, then the essential norm 
of the Toeplitz operator $T(a)$ with $a\in C+H^\infty$ on the weighted 
Hardy space $H^2(w)$ is equal to $\|a\|_{L^\infty}$.
\end{abstract}
\subjclass{Primary 47B35, 46E30}
\keywords{Toeplitz operator, essential norm, weighted Hardy space,
essential norm, $C+H^\infty$.}
\maketitle
\section{Introduction}
For Banach spaces $\mathcal{X},\mathcal{Y}$, let 
$\mathcal{B}(\mathcal{X},\mathcal{Y})$ and 
$\mathcal{K}(\mathcal{X},\mathcal{Y})$ 
denote the sets of bounded linear and compact linear operators from 
$\mathcal{X}$ to $\mathcal{Y}$, respectively. The norm of an operator 
$A\in\mathcal{B}(\mathcal{X},\mathcal{Y})$ is denoted
by $\|A\|_{\mathcal{B}(\mathcal{X},\mathcal{Y})}$. The essential norm of 
$A \in \mathcal{B}(\mathcal{X},\mathcal{Y})$ 
is defined as follows:
\[
\|A\|_{\mathcal{B}(\mathcal{X},\mathcal{Y}),\mathrm{e}} 
:= 
\inf\{\|A - K\|_{\mathcal{B}(\mathcal{X},\mathcal{Y})}\ : \  K \in 
\mathcal{K}(\mathcal{X},\mathcal{Y})\}.
\]
As usual, we abbreviate $\mathcal{B}(\mathcal{X},\mathcal{X})$ and 
$\mathcal{K}(\mathcal{X},\mathcal{X})$ to
$\mathcal{B}(\mathcal{X})$ and $\mathcal{K}(\mathcal{X})$, respectively.

Let $\mathbb{T}:=\{z\in\mathbb{C}:|z|=1\}$ be the unit circle in the complex
plane. We equip $\mathbb{T}$ with  the Lebesgue measure $m$ normalised so 
that $m(\mathbb{T})=1$. In this paper, all function spaces will be 
considered over $\mathbb{T}$.
A measurable function $w:\mathbb{T}\to[0,\infty]$ is said to be a weight
if $0<w<\infty$ a.e. on $\mathbb{T}$. Let $1<p<\infty$ and $w$ be a weight.
Weighted Lebesgue spaces $L^p(w)$ consist of all measurable functions
$f:\mathbb{T}\to\mathbb{C}$ such that $fw\in L^p$. The norm in $L^p(w)$
is defined by
\[
\|f\|_{L^p(w)}:=\|fw\|_{L^p}
=
\left(\int_{\mathbb{T}}|f(t)|^pw^p(t)\,dm(t)\right)^{1/p}.
\]
Consider the operators $S$ and $P$, defined for a function 
$f\in L^1$ and a.e. point $t\in\mathbb{T}$ by
\[
(Sf)(t):=\frac{1}{\pi i}\,\mbox{p.v.}\int_\mathbb{T}
\frac{f(\tau)}{\tau-t}\,d\tau,
\quad
(Pf)(t):=\frac{1}{2}(f(t)+(Sf)(t)),
\]
respectively, where the integral is understood in the Cauchy principal 
value sense. The operator $S$ is called the Cauchy singular 
integral operator and the operator $P$ is called the Riesz projection.
It is well known (see \cite[Theorem~1]{HMW73}
and also \cite[Thorem~4.15]{BK97}, \cite[Ch.~VI]{G07}) that the 
Riesz projection $P$ is bounded on $L^p(w)$ if and only if $w$ 
belongs to the Muckenhoupt class $A_p$, that is,
\[
\sup_{\gamma\subset\mathbb{T}}
\left(\frac{1}{m(\gamma)}\int_\gamma w^p(t)\, dm(t)\right)^{1/p}
\left(\frac{1}{m(\gamma)}\int_\gamma w^{-p'}(t)\, dm(t)\right)^{1/p'}
<\infty,
\]
where $1/p+1/p'=1$ and the supremum is taken over all arcs 
$\gamma\subset\mathbb{T}$. If $w\in A_p$, then 
$L^\infty\hookrightarrow L^p(w)\hookrightarrow L^1$.

For a function $f \in L^1$, let
\[
\widehat{f}(n) 
= 
\frac{1}{2\pi} \int_{-\pi}^\pi 
f\left(e^{i\theta}\right) e^{-i n\theta}\, d\theta , 
\quad 
n \in \mathbb{Z}
\]
be the Fourier coefficients of $f$. For $1<p<\infty$ and $w\in A_p$, let
\[
H^p(w) = H[L^p(w)] := 
\left\{
f\in L^p(w)\ :\ \widehat{f}(n)=0\quad\mbox{for all}\quad n<0
\right\} 
\]
be the weighted Hardy space. The classical Hardy spaces  $H^p$, 
$1\le p\le\infty$, are defined similarly if one replaces $L^p(w)$ by $L^p$ 
in the above definition. It is well known that if $w\in A_p$, then $P$ maps 
$L^p(w)$ onto $H^p(w)$.

Let $1<p<\infty$ and $w\in A_p$. For $a\in L^\infty$, the Toeplitz operator 
with symbol $a$ is defined by
\[
T(a)f=P(af),
\quad 
f\in H^p(w).
\]
It is clear that $T(a)\in\mathcal{B}(H^p(w))$ and 
\[
\|T(a)\|_{\mathcal{B}(H^p(w)),{\rm e}}
\le 
\|T(a)\|_{\mathcal{B}(H^p(w))}
\le 
\|P\|_{\mathcal{B}(L^p(w))}\|a\|_{L^\infty}.
\]

I.~Gohberg and N.~Krupnik \cite[Theorem~6]{GK68} proved that 
$\|P\|_{\mathcal{B}(L^p),\mathrm{e}}\ge 1/\sin(\pi/p)$ and conjectured that
$\|P\|_{\mathcal{B}(L^p)}=1/\sin(\pi/p)$. This conjecture was confirmed
by B.~Hollenbeck and I.~Verbitsky in \cite{HV00}. Thus
\begin{equation}\label{eq:Toeplitz-bound-Hp}
\|T(a)\|_{\mathcal{B}(H^p),\mathrm{e}}
\le 
1/\sin(\pi/p)
\|a\|_{L^\infty},
\quad
a\in L^\infty.
\end{equation}

Let $C$ denote the Banach space of all complex-valued continuous functions
on $\mathbb{T}$ with the supremum norm and let 
\[
C+H^\infty:=\{f\in L^\infty\ :\ f=g+h,\ g\in C,\ h\in H^\infty\}.
\]
In 1967, Sarason observed that $C+H^\infty$ is a closed subalgebra
of $L^\infty$ (see, e.g., \cite[Ch.~IX, Theorem~2.2]{G07} for the proof 
of this fact).

Consider power weights of the form
\begin{equation}\label{eq:Khvedelidze}
\varrho(t):=\prod_{j=1}^k |t-t_j|^{\lambda_j}
\end{equation}
where $t_1,\dots,t_k\in\mathbb{T}$ are pairwise disjoint and 
$\lambda_1,\dots,\lambda_k\in\mathbb{R}$. These weights are usually called
Khvedelidze weights. It is well known that $\varrho\in A_p$ with
$1<p<\infty$ if and only if $-1/p<\lambda_j<1-1/p$ for all $j\in\{1,\dots,k\}$
(see, e.g., \cite[Theorem~2.2]{BK97}).

In 1988, A.~B\"ottcher, N.~Krupnik, and B.~Silbermann proved the following
(see \cite[Corollary~7.5]{BKS88}).
\begin{theorem}\label{th:BKS}
If $\varrho$ is a power weight of the form \eqref{eq:Khvedelidze} satisfying
\[
-1/2<\lambda_j<1/2 
\quad\mbox{for all}\quad
j\in\{1,\dots,k\}, 
\]
then for all $a\in C$
one has $\|T(a)\|_{\mathcal{B}(H^2(\varrho)),\mathrm{e}}=\|a\|_{L^\infty}$.
\end{theorem}
In particular, the essential norm of $T(a)$ on $H^2(\varrho)$ does not
depend on the choice of a Khvedelidze weight
$\varrho\in A_2$. Further they asked in \cite[Section~7.6]{BKS88} 
whether the essential norm of Toeplitz operators $T(a)$ with 
$a\in C$ acting on Hardy spaces $H^p$ depend on $p\in(1,\infty)$. 
The second author answered this question in the negative \cite{S21}. 
More precisely, it was shown that the equality 
$\|T(a)\|_{\mathcal{B}(H^p),\mathrm{e}}=\|a\|_{L^\infty}$ 
holds for all $a\in C$ if and only if $p=2$. Nevertheless, the following 
estimates for $\|T(a)\|_{\mathcal{B}(H^p),\mathrm{e}}$ hold
(see \eqref{eq:Toeplitz-bound-Hp} and \cite[Theorem~4.1]{S21}).
\begin{theorem}\label{th:Shargorodsky-JFA}
Let $1<p<\infty$ and $a\in C+H^\infty$. Then
\[
\|a\|_{L^\infty} 
\le 
\|T(a)\|_{\mathcal{B}(H^p),\mathrm{e}}  
\le 
\min\left\{2^{|1-2/p|},1/\sin(\pi/p)\right\} \|a\|_{L^\infty}.
\]
\end{theorem}

The aim of this paper it to study the relations between the essential norms of 
Toeplitz operators with symbols in $C+H^\infty$on unweighted and weighted 
Hardy spaces. The following is our main result.
\begin{theorem}\label{th:main}
Let $1<p<\infty$ and $w\in A_p$. If $a\in C+H^\infty$, then
\[
\|T(a)\|_{\mathcal{B}(H^p),\mathrm{e}}
=
\|T(a)\|_{\mathcal{B}(H^p(w)),\mathrm{e}}.
\]
\end{theorem}
It is instructive to compare this result with the behaviour of the essential 
norm $\|S\|_{\mathcal{B}(L^p(\varrho)),\mathrm{e}}$ of the Cauchy singular 
integral operator with $\varrho\in A_p$ of the form \eqref{eq:Khvedelidze}, 
which depends not only on $p\in (1,\infty)$, but also on $\varrho$ 
(see \cite[Remark~6.8]{BKS88} and Krupnik's survey \cite{K10}).

In fact, we prove an abstract form of Theorem~\ref{th:main} 
(see Theorem~\ref{th:independence-of-weight}), where we replace the Hardy 
space $H^p$ with $1<p<\infty$ by the abstract Hardy space $H[X]$ built
upon a Banach function space $X$ such that $P\in\mathcal{B}(X)$
and replace the weighted Hardy space $H^p(w)$ with $w\in A_p$ by
the abstract Hardy space $H[X(w)]$ built upon 
$X(w)=\{f:\mathbb{T}\to\mathbb{C}:fw\in X\}$ with a weight $w$ such that
$P\in\mathcal{B}(X(w))$. Our proof is fairly elementary and is completely 
different from that of Theorem~\ref{th:BKS}.

It is clear that Theorems~\ref{th:Shargorodsky-JFA} and~\ref{th:main}
immediately imply the following extension of Theorem~\ref{th:BKS}.
\begin{corollary}
Suppose that $1<p<\infty$ and $w\in A_p$. If $a\in C+H^\infty$, then
\[
\|a\|_{L^\infty} 
\le 
\|T(a)\|_{\mathcal{B}(H^p(w)),\mathrm{e}}  
\le 
\left\{2^{|1-2/p|},1/\sin(\pi/p)\right\} \|a\|_{L^\infty}.
\]
In particular, if $w\in A_2$ and $a\in C+H^\infty$, then 
$\|T(a)\|_{\mathcal{B}(H^2(w)),\mathrm{e}}=\|a\|_{L^\infty}$.
\end{corollary}

We will use the following notation:
\[
\mathbf{e}_m(z) := z^m, 
\quad z \in \mathbb{C}, 
\quad m \in \mathbb{Z}.
\]

The paper is organised as follows. In Section~\ref{sec:preliminaries}, we 
collect definitions of a Banach function space, its associate space $X'$, and 
a weighted Banach function space $X(w)$. Note that if $w\in X$ and $1/w\in X'$, 
then $X(w)$ is a Banach function space itself. This allows us to define
the abstract Hardy space $H[X(w)]$ built upon $X(w)$. We recall that if
$w\in X$ and $1/w\in X'$, then then the abstract Hardy spaces $H[X]$ and
$H[X(w)]$ are isometrically isomorphic. We conclude the preliminaries with
some useful properties of the Riesz projection $P$ and the definition of
a Toeplitz operator $T(a)$ on the abstract Hardy space built upon a Banach
function space $X$ such that $P\in\mathcal{B}(X)$.

Section~\ref{sec:proof} is devoted to the proof of Theorem~\ref{th:main}.
First we recall an observation from \cite{OKES20} that the Toeplitz operator
$T(\mathbf{e}_{-n}h)$ with $n\in\mathbb{N}$ and $h\in H^\infty$
is bounded on the abstract Hardy space $H[X]$ for an arbitrary Banach function
space, i.e. even without the assumption that $P\in\mathcal{B}(X)$.
Further, we show that the essential norms of $T(\mathbf{e}_{-n}h)$
on $H[X]$ and $H[X(w)]$ coincide if $w\in X$ and $1/w\in X'$.
Here we essentially use the isomorphic isomorphism of $H[X]$ and $H[X(w)]$.
Finally, taking into account that the set 
$\{\mathbf{e}_{-n}h:n\in\mathbb{N},h\in H^\infty\}$ is dense in $C+H^\infty$, 
we prove Theorem~\ref{th:independence-of-weight}. In turn, 
Theorem~\ref{th:independence-of-weight} yield Theorem~\ref{th:main}.
\section{Preliminaries}\label{sec:preliminaries}
\subsection{Banach function spaces}\label{subsec:BFS}
Let $\mathcal{M}$ be the set of all measurable complex-valued functions on 
$\mathbb{T}$ equipped with the normalized Lebesgue measure $m$ and let 
$\mathcal{M}^+$ be the subset of functions in $\mathcal{M}$ whose values lie 
in $[0,\infty]$. Following \cite[Ch.~1, Definition~1.1]{BS88}, a mapping 
$\rho: \mathcal{M}^+\to [0,\infty]$ is called a Banach function norm if, 
for all functions $f,g, f_n\in \mathcal{M}^+$ with $n\in\mathbb{N}$, and for 
all constants $a\ge 0$, the following  properties hold:
\begin{eqnarray*}
{\rm (A1)} & &
\rho(f)=0  \Leftrightarrow  f=0\ \mbox{a.e.},
\
\rho(af)=a\rho(f),
\
\rho(f+g) \le \rho(f)+\rho(g),\\
{\rm (A2)} & &0\le g \le f \ \mbox{a.e.} \ \Rightarrow \ 
\rho(g) \le \rho(f)
\quad\mbox{(the lattice property)},\\
{\rm (A3)} & &0\le f_n \uparrow f \ \mbox{a.e.} \ \Rightarrow \
       \rho(f_n) \uparrow \rho(f)\quad\mbox{(the Fatou property)},\\
{\rm (A4)} & & \rho(1) <\infty,\\
{\rm (A5)} & &\int_{\mathbb{T}} f(t)\,dm(t) \le C\rho(f)
\end{eqnarray*}
with {a constant} $C \in (0,\infty)$ that may depend on $\rho$,  but is 
independent of $f$. When functions differing only on a set of measure zero 
are identified, the set $X$ of all functions $f\in \mathcal{M}$ for which 
$\rho(|f|)<\infty$ is called a Banach function space. For each $f\in X$, 
the norm of $f$ is defined by $\|f\|_X :=\rho(|f|)$. The set $X$ equipped 
with the natural vector space operations and this norm becomes a Banach 
space (see \cite[Ch.~1, Theorems~1.4 and~1.6]{BS88}). If $\rho$ is a Banach 
function norm, its associate norm $\rho'$ is defined on $\mathcal{M}^+$ by
\[
\rho'(g):=\sup\left\{
\int_\mathbb{T} f(t)g(t)\,dm(t) \ : \ 
f\in\mathcal{M}^+, \ \rho(f) \le 1
\right\}, \quad g\in \mathcal{M}^+.
\]
It is a Banach function norm itself \cite[Ch.~1, Theorem~2.2]{BS88}.
The Banach function space $X'$ determined by the Banach function 
norm $\rho'$ is called the associate space (K\"othe dual) of $X$. 
The associate space $X'$ can be viewed as a subspace of the 
Banach dual space $X^*$. 
\subsection{Weighted Banach function spaces}
For a weight $w$ and a Banach function space $X$, the weighted space $X(w)$ 
consists of all measurable functions $f:\mathbb{T}\to\mathbb{C}$ such that 
$fw\in X$. We equip it with the norm $\|f\|_{X(w)}=\|fw\|_X$.
\begin{lemma}\label{le:weighted-BFS}
Let $X$ be a Banach function space with the associate spaces $X'$ and 
let $w$ be a weight.

\begin{enumerate}
\item[(a)]
If $w\in X$ and $1/w\in X'$, then $X(w)$ is a Banach function space,
whose associate space is $X'(1/w)$.

\item[(b)]
If $P$ is bounded on the normed space $X(w)$, then $w\in X$ and 
$1/w\in X'$.
\end{enumerate}
\end{lemma}
Part (a) was proved in \cite[Lemma~2.3]{K03}, part (b) follows from 
\cite[Theorem~6.1]{K03}. 
\subsection{Isometric isomorphism of weighted and nonweighted abstract
Hardy spaces}
Let $X$ be a Banach function space and let $X'$ be its associate space.
Lem\-ma~\ref{le:weighted-BFS}(a) implies that if $w$ is a weight such that
$w\in X$ and $1/w\in X'$, then $X(w)$ is a Banach function space. Then 
$L^\infty\hookrightarrow X(w)\hookrightarrow L^1$
and the abstract Hardy space $H[X(w)]$ built upon $X(w)$ is well defined
by
\[
H[X(w)]:=
\big\{
f\in X(w)\ :\ \widehat{f}(n)=0\quad\mbox{for all}\quad n<0
\big\}. 
\]
For simplicity, we will write $H[X]$ if $w=1$.
\begin{lemma}[{\cite[Lemma~4.2]{KS24}}]
\label{le:isometric-isomorphism}
Let $X$ be a Banach function space with the associate space $X'$
and let $w$ be a weight such that $w\in X$ and $1/w\in X'$.
Consider
\[
W(z) := \exp\left(\frac1{2\pi}\int_{-\pi}^\pi
\frac{e^{it} + z}{e^{it} - z}\, \log w(e^{it})\, dt\right), 
\quad 
|z|<1.
\]
Then the mapping $f\mapsto M_Wf:=W\cdot f$
is an isometric isomorphism of $H[X(w)]$ onto $H[X]$ and
the mapping $g\mapsto M_{W^{-1}}g :=W^{-1}\cdot g$ is an isometric
isomorphism of $H[X]$ onto $H[X(w)]$.
\end{lemma}
\subsection{Auxiliary lemmas on the Riesz projection}
We will need the following auxiliary lemmas.
\begin{lemma}[{\cite[Lemma~3.1]{K21}}]
\label{le:Pf=g}
Let $f \in L^1$. Suppose there exists $g \in H^1$ such that 
$\widehat{f}(n) = \widehat{g}(n)$ for all $n \ge 0$. Then $Pf = g$.
\end{lemma}
The following lemma justifies the terminology about the operator $P$.
\begin{lemma}[{\cite[Lemma~1.1]{KS19-JMAA}}]
\label{le:Riesz-projection-properties}
If $X$ is a Banach function space on which the operator $P$ is bounded,
then $P$ maps the space $X$ onto the abstract Hardy space $H[X]$.
\end{lemma}
\subsection{Toeplitz operators on abstract Hardy spaces}
Let $X$ be a Banach function space. If $a\in L^\infty$, then the multiplication
operator by $a$ defined by
\[
M_af:=af, \quad f\in X,
\]
is bounded on $X$ and $\|M_a\|_{\mathcal{B}(X)}\le \|a\|_{L^\infty}$.
Further, if the Riesz projection is bounded on $X$, then the Toeplitz
operator defined by
\[
T(a)f:=P(M_af)=P(af),
\quad
f\in H[X],
\]
is bounded on $H[X]$ and
\begin{equation}\label{eq:norm-Toeplitz}
\|T(a)\|_{\mathcal{B}(H[X])}\le \|P\|_{\mathcal{B}(X)}\|a\|_{L^\infty}.
\end{equation}
\section{Proof of the main result}\label{sec:proof}
\subsection{Operator $P_n$ and a special Toeplitz operator}
We will need auxiliary results on a representation of a Toeplitz
operator $T(\mathbf{e}_{-n}h)$ with $n\in\mathbb{N}$ and $h\in H^\infty$
in terms of the projection $P_n$ from $X$ onto the subspace of $H[X]$ of 
analytic polynomials of order $n-1$ obtained in our recent paper \cite{OKES20}.
We would like to underline that we do not require here and in
Subsection~\ref{subsec:independce-of-weight-for-special-symbols} 
that the Riesz projection $P$ is bounded on $X$.
\begin{lemma}[{\cite[Lemma~3.1, Corollary~3.2]{OKES20}}]
\label{le:boundedness-Pn}
For every $n\in\mathbb{N}$ and $f\in L^1$, let
\[
P_n f := \sum_{k = 0}^{n-1} \widehat{f}(k) \mathbf{e}_k.
\]
Then the operator $P_n$ is bounded from $L^1$ to $H^\infty$ and
from a Banach function space $X$ to the abstract Hardy space $H[X]$.
\end{lemma}
\begin{lemma}[{\cite[Lemma~3.4 and equality (3.4)]{OKES20}}]
\label{le:boundedness-T-special}
Let $X$ be a Banach function space. If $n \in \mathbb{Z}_+$ and $h \in H^\infty$, 
then the Toeplitz operator $T(\mathbf{e}_{-n} h) : H[X] \to H[X]$ is bounded
and
\begin{equation}\label{eq:boundendess-T-special}
T(\mathbf{e}_{-n} h) f 
= 
P(\mathbf{e}_{-n} h f) = \mathbf{e}_{-n} (I - P_n)(hf),
\quad
f\in H[X].
\end{equation} 
\end{lemma}
\subsection{The essential norm of a special Toeplitz operator on a weighted 
Hardy space is independent of the weight}
\label{subsec:independce-of-weight-for-special-symbols}
In this subsection we show that the essential norms of the special
Toepllitz operators with symbols of the form 
$\mathbf{e}_{-n}h$, where $n\in\mathbb{N}$ and $h\in H^\infty$,
are independent of the weight.
\begin{lemma}\label{le:special-T-independence-of-weight} 
Let $X$ be a Banach function space with the associate space $X'$
and let $w$ be a weight satisfying $w\in X$ and $1/w\in X'$. 
If $n \in \mathbb{N}$ and $h \in H^\infty$, then
\[
\|T(\mathbf{e}_{-n} h)\|_{\mathcal{B}(H[X(w)]),\mathrm{e}} 
= 
\|T(\mathbf{e}_{-n} h)\|_{\mathcal{B}(H[X]),\mathrm{e}} .
\]
\end{lemma}
\begin{proof}
Let the operators $M_W$ and $M_{W^{-1}}$ be as in 
Lemma~\ref{le:isometric-isomorphism}. Since $W$ is an outer function 
(see \cite[Ch.~5]{H88}), one has $|W|=w$ 
a.e. on $\mathbb{T}$, and $W\in H[X]$, $W^{-1}\in H[X']$.
 
By Lemma~\ref{le:boundedness-T-special},
the operator $T(\mathbf{e}_{-n}h)$ is bounded on $H[X]$ and $H[X(w)]$.
For any $f \in H[X]$, one gets from \eqref{eq:boundendess-T-special}
and Lemma~\ref{le:Pf=g} that
\begin{align*}
&  M_{W} T(\mathbf{e}_{-n} h) M_{W^{-1}} f 
= 
W \mathbf{e}_{-n} (I - P_n)\left(hW^{-1} f\right) 
\\
& \quad = 
P\left(W \mathbf{e}_{-n} (I - P_n)\left(hW^{-1} f\right)\right) 
\\
& \quad = 
P\left(W \mathbf{e}_{-n} (hW^{-1} f) 
- 
W \mathbf{e}_{-n} P_n\left(hW^{-1} f\right)\right) 
\\
& \quad = 
P\left(\mathbf{e}_{-n} (h f) 
- 
W \mathbf{e}_{-n} P_n \left(hW^{-1} f\right)\right) 
\\
& \quad = 
P\left(\mathbf{e}_{-n}  (I - P_n)(h f) 
+ 
\mathbf{e}_{-n} P_n(h f) 
- 
W \mathbf{e}_{-n} P_n\left(hW^{-1} f\right)\right) 
\\
& \quad = 
P\left(T(\mathbf{e}_{-n} h) f 
+  
\mathbf{e}_{-n} P_n (h f) 
- 
W \mathbf{e}_{-n} P_n \left(hW^{-1} f\right)\right) 
\\
& \quad = 
T(\mathbf{e}_{-n} h) f 
+ 
T(\mathbf{e}_{-n}) P_n  M_h f 
- 
T(\mathbf{e}_{-n}) M_W P_n M_{hW^{-1}} f .
\end{align*}

Since $W^{-1} \in X'$ and $h\in H^\infty$, it follows from the H\"older 
inequality for Banach function spaces (see \cite[Ch.~1, Theorem~2.4]{BS88})
that $M_{hW^{-1}} \in \mathcal{B}(H[X], L^1)$. 

On the other hand, taking into account 
\cite[Section~3.3.1, properties (a), (g)]{N19}, since 
$W\in H[X]\subset H^1$, one has $Wg\in H^1\cap X=H[X]$
for every $g\in H^\infty$.
Hence $M_W \in \mathcal{B}(H^\infty, H[X])$. 

Then Lemmas~\ref{le:boundedness-Pn} and  \ref{le:boundedness-T-special} 
(with $h = 1$) imply that
\[
T(\mathbf{e}_{-n}) M_W P_n M_{hW^{-1}} \in \mathcal{B}(H[X]),
\quad
T(\mathbf{e}_{-n}) P_n M_h \in \mathcal{B}(H[X]) .
\]
So,
\[
K_0  :=  
T(\mathbf{e}_{-n}) P_n M_h 
- 
T(\mathbf{e}_{-n}) M_W P_n M_{hW^{-1}} 
\]
is a bounded finite-rank operator on $H[X]$, and
\begin{equation}\label{K0}
M_{W} T(\mathbf{e}_{-n} h) M_{W^{-1}} = T(\mathbf{e}_{-n} h) + K_0 .
\end{equation}
Since $K_0 \in \mathcal{K}(H[X])$, one has
\[
M_{W} K M_{W^{-1}} - K_0 \in \mathcal{K}(H[X]) 
\]
for every $K \in \mathcal{K}(H[X(w)])$. Moreover, for every 
$T \in \mathcal{K}(H[X])$, there exists $K \in \mathcal{K}(H[X(w)])$ such 
that
\[
T = M_{W} K M_{W^{-1}} - K_0 .
\]
Indeed, this $K$ is given by the formula
\[
K = M_{W^{-1}}(T + K_0) M_{W} \in \mathcal{K}(H[X(w)]) .
\]
Using \eqref{K0}, one gets
\begin{align*}
&
\|T(\mathbf{e}_{-n} h)\|_{\mathcal{B}(H[X(w)]),\mathrm{e}} 
= 
\inf_{K \in \mathcal{K}(H[X(w)])} 
\|T(\mathbf{e}_{-n} h) - K\|_{\mathcal{B}(H[X(w)])}  
\\
&\quad = 
\inf_{K \in \mathcal{K}(H[X(w)])} 
\|M_{W}(T(\mathbf{e}_{-n} h) - K)M_{W^{-1}}\|_{\mathcal{B}(H[X])}  
\\
&\quad = 
\inf_{K \in \mathcal{K}(H[X(w)])} 
\|T(\mathbf{e}_{-n} h) - (M_{W} K M_{W^{-1}} - K_0)\|_{\mathcal{B}(H[X])}  
\\
&\quad = 
\inf_{T \in \mathcal{K}(H[X])} 
\|T(\mathbf{e}_{-n} h) - T\|_{\mathcal{B}(H[X])} 
= 
\|T(\mathbf{e}_{-n} h)\|_{\mathcal{B}(H[X]),\mathrm{e}}, 
\end{align*}
which completes the proof.
\end{proof}
\subsection{The essential norm of a Toeplitz operator on a weighted
Hardy space is independent of the weight}
Now we are in a position to prove the main result of this section.
\begin{theorem}\label{th:independence-of-weight}
Let $X$ be a Banach function space and let $w$ be a weight such that 
the Riesz projection $P$ is bounded on the spaces $X$ and $X(w)$.
If $a\in C+H^\infty$, then 
\begin{equation}\label{eq:independence-of-weight-1}
\|T(a)\|_{\mathcal{B}(H[X(w)]),\mathrm{e}} 
= 
\|T(a)\|_{\mathcal{B}(H[X]),\mathrm{e}}.
\end{equation}
\end{theorem}
\begin{proof}
It follows from \cite[Ch.~IX, Theorem~2.2]{G07} that there is a sequence
$\{a_m\}$ of functions of the form $\mathbf{e}_{-n}h$ with $h \in H^\infty$
and $n \in \mathbb{N}$ such that $\|a-a_m\|_{L^\infty}\to 0$ as $m\to\infty$.
Since $P\in\mathcal{B}(X(w))$, Lemma~\ref{le:weighted-BFS}(b) implies that
$w\in X$ and $1/w\in X'$. Then Lemma~\ref{le:special-T-independence-of-weight} 
yields 
\begin{equation}\label{eq:independence-of-weight-2}
\|T(a_m)\|_{\mathcal{B}(H[X(w)]),\mathrm{e}} 
= 
\|T(a_m)\|_{\mathcal{B}(H[X]),\mathrm{e}},
\quad
m\in\mathbb{N}.
\end{equation}
Since $P\in\mathcal{B}(X)$, we obtain from \eqref{eq:norm-Toeplitz}
that
\begin{align*}
\left|
\|T(a)\|_{\mathcal{B}(H[X]),\mathrm{e}}
-
\|T(a_m)\|_{\mathcal{B}(H[X]),\mathrm{e}}
\right|
& \le 
\|T(a-a_m)\|_{\mathcal{B}(H[X]),\mathrm{e}}
\\
& \le 
\|P\|_{\mathcal{B}(X)}\|a-a_m\|_{L^\infty},
\end{align*}
whence
\begin{equation}\label{eq:independence-of-weight-3}
\lim_{m\to\infty}\|T(a_m)\|_{\mathcal{B}(H[X]),\mathrm{e}}
=
\|T(a)\|_{\mathcal{B}(H[X]),\mathrm{e}}.
\end{equation}
Similarly, since $P\in\mathcal{B}(X(w))$, we get
\begin{equation}\label{eq:independence-of-weight-4}
\lim_{m\to\infty}\|T(a_m)\|_{\mathcal{B}(H[X(w)]),\mathrm{e}}
=
\|T(a)\|_{\mathcal{B}(H[X(w)]),\mathrm{e}}.
\end{equation}
Combining 
\eqref{eq:independence-of-weight-2}--\eqref{eq:independence-of-weight-4},
we arrive at \eqref{eq:independence-of-weight-1}.
\end{proof}
\subsection{Proof of Theorem~\ref{th:main}}
If $1<p<\infty$ and $w\in A_p$, then the Riesz projection $P$ is bounded
on the Lebesgue space $L^p$ and on its weighted counterpart $L^p(w)$.
In this case, Theorem~\ref{th:main} follows immediately from 
Theorem~\ref{th:independence-of-weight}.
\qed
\section*{Declarations}
\subsection*{Funding}
This work is funded by national funds through the FCT - Funda\c{c}\~ao para a 
Ci\^encia e a Tecnologia, I.P., under the scope of the projects UIDB/00297/2020 
(\url{https://doi.org/10.54499/UIDB/00297/2020})
and UIDP/ 00297/2020

\noindent 
(\url{https://doi.org/10.54499/UIDP/00297/2020})
(Center for Mathematics and Applications).
\subsection*{Data availability}
Data sharing not applicable to this article as no datasets
were generated or analysed during the current study.
\subsection*{Conflict of interest}
All authors certify that they have no affiliations with
or involvement in any organisation or entity with any financial interest or
non-financial interest in the subject matter or materials discussed in this
manuscript.
\bibliographystyle{abbrv}
\bibliography{OKES21}
\end{document}